\documentclass[10pt,a4,fleqn]{article}
\usepackage{amsmath, amsthm, amssymb}
\usepackage{graphicx}

\newtheorem{thm}{Theorem}
\newtheorem{cor}[thm]{Corollary}
\newtheorem{lem}[thm]{Lemma}
\newtheorem{prop}[thm]{Proposition}
\newtheorem{question}{Question}
\theoremstyle{definition}
\newtheorem{defn}{Definition}
\theoremstyle{remark}
\newtheorem{rem}[thm]{Remark}

\newcommand{\comment}[1]{}
\newif \ifskip
\skiptrue

\begin{document}
\begin{titlepage}
\title{A Most General Edge Elimination Graph Polynomial}

\author{
Ilia Averbouch\thanks{ Partially supported by a grant of the
Graduate School of the Technion--Israel Institute of Technology }
\mbox{\ \ \ }
Benny Godlin\thanks{ Partially supported by a grant of the Graduate
School of the Technion--Israel Institute of Technology }
\mbox{\ \ \ }
J.A.Makowsky\thanks{ Partially supported by a grant of the Fund for
Promotion of Research of the Technion--Israel Institute of
Technology and a grant of the Israel Science Foundation (2007-2010)}
\\
\ \\
Faculty of Computer Science\\
Israel Institute of Technology\\
Haifa, Israel
\\
\{ailia,bgodlin,janos\}@cs.technion.ac.il }

\maketitle

\begin{abstract}
We look for graph polynomials which satisfy recurrence relations on
three kinds of edge elimination:  edge deletion, edge contraction
and deletion  of edges together with their end points. Like in the
case of deletion and contraction only (W. Tutte, 1954), it turns out
that there is a most general polynomial satisfying such recurrence
relations, which we call $\xi(G,x,y,z)$. We show that the new
polynomial simultaneously generalizes the Tutte polynomial, the
matching polynomial, and the recent generalization of the chromatic
polynomial proposed by K.Dohmen, A.P\"{o}nitz and P.Tittman (2003),
including also the independent set polynomial of 
I. Gutman and F. Harary, (1983) and
the vertex-cover polynomial of 
F,M. Dong, M.D. Hendy, K.T. Teo and C.H.C. Little (2002).
We establish two definitions of the new polynomial:
first, the most general confluent recursive definition, and then an
explicit one, using a set expansion formula, and prove their identity.
We further expand this result to edge-labeled graphs as was done
for the Tutte polynomial by T. Zaslavsky (1992) and  
B. Bollob\'as and O. Riordan (1999).
The edge labeled polynomial $\xi_{lab}(G,x,y,z, \bar{t})$
also generalizes the chain polynomial of 
R.C. Read and E.G. Whitehead Jr. (1999).
Finally, we discuss the complexity of computing $\xi(G,x,y,z)$.
\end{abstract}
\end{titlepage}
\maketitle
\section{Introduction}
There are several well-studied graph polynomials, among them the
chromatic polynomial,
\cite{bk:Biggs93,bk:GodsilRoyle01,bk:DongKohTeo2005}, different
versions of the Tutte polynomial,
\cite{bk:Bollobas99,ar:BollobasRiordan99,ar:Sokal2005a}, and of the
matching polynomial,
\cite{ar:HeilmannLieb72,bk:LovaszPlummer86,bk:GodsilRoyle01}, which
are known to satisfy certain linear recurrence relations with
respect to \emph{deletion} of an edge, \emph{contraction} of an
edge, or deletion of an edge together with its endpoints, which we
call \emph{extraction} of an edge. The generalization of the
chromatic polynomial, which was introduced by K.Dohmen, A.P\"{o}nitz
and P.Tittman in \cite{ar:DPT03}, happens to satisfy such recurrence
relation as well. The question that arises is, what is the most
general graph polynomial that satisfies similar linear recurrence
relation.

In this paper all the graphs are unlabeled unless it is explicitly
mentioned; multiple edges and self loops are allowed. We denote by
$G=(V,E)$ the graph with vertex set $V$ and edge set $E$.

\subsection{Recursive definition of graph polynomials}

We define three basic edge elimination operations on multigraphs:
\begin{itemize}
\item \emph{Deletion}. We denote by $G_{-e}$ the graph obtained from $G$ by
simply removing the edge $e$.
\item \emph{Contraction}. We denote by $G_{/e}$ the graph obtained
from $G$ by unifying the endpoints of $e$. Note that this operation
can cause production of multiple edges and self loops.
\item \emph{Extraction}. We denote by $G_{\dagger e}$ the graph induced by
$V\setminus\{u,v\}$ provided $e=\{u,v\}$. Note that this operation
removes also all the edges adjacent to $e$.
\end{itemize}
Additionally,  we require the polynomial to be
\emph{multiplicative}
for disjoint unions,
i.e., if $G_1 \oplus G_2$ denotes disjoint
union of two graphs, then the polynomial $P(G_1 \oplus G_2) = P(G_1)
\cdot P(G_2)$.
This is justified by the fact that
the polynomials occurring in the literature are usually
multiplicative.
The initial conditions are defined for an \emph{empty
set} (graph without vertices, usually, $P(\emptyset)=1$) and for a
single point $P(E_1)$. With respect to these operations, we recall
the known recursive definitions of graph polynomials:

\paragraph{Matching polynomial.}
There are different versions of the matching polynomial discussed in
the literature, for example \emph{matching generating polynomial}
$g(G,\lambda)=\sum_{i=0}^n a_i \lambda^i$ and \emph{matching defect
polynomial} $\mu(G,\lambda)=\sum_{i=0}^n (-1)^i a_i \lambda^{n-2i}$,
where $n=|V|$ and $a_i$ is the number of $i$-matchings in $G$. We
shall use the bivariate version that incorporates the both above:
\begin{equation} \label{general_matching_eqn}
M(G,x,y)=\sum_{i=0}^n a_i x^{n-2i} y^i
\end{equation}
The recursive definition of this polynomial is as follows:
\begin{eqnarray}
\nonumber M(G) &=& M(G_{-e}) + y \cdot M(G_{\dagger e})\\
\nonumber M(G_1 \oplus G_2) &=& M(G_1)\cdot M(G_2)\\
\nonumber M(E_1)&=&x;\\
M(\emptyset)&=&1; \label{rec_matching}
\end{eqnarray}

\paragraph{Tutte polynomial.}
We recall the definition of classical two-variable Tutte polynomial
(cf. for example B.Bollob\'as \cite{bk:Bollobas99}):
\begin{defn}
\emph{Let $G=(V,E)$ be a (multi-)graph. Let $A\subseteq E$ be a
subset of edges. We denote by $k(A)$ the number of connected
components in the spanning subgraph $(V,A)$. Then two-variable Tutte
polynomial is defined as follows}
\end{defn}
\begin{equation}\label{tutte_cls}
T(G,x,y)=\sum_{A \subseteq E}(x-1)^{k(A)-k(E)}(y-1)^{|A|+k(A)-|V|}
\end{equation}
This polynomial has linear recurrence relation with respect to the
operations above:
\begin{eqnarray}
\nonumber T(G,x,y)&=&\left\{\begin{array}{lll}
                x \cdot T(G_{/e},x,y) & if~e~is~a~bridge, \\
                y \cdot T(G_{-e},x,y) & if~e~is~a~loop, \\
                T(G_{/e},x,y)+T(G_{-e},x,y) & otherwise
                \end{array}
                \right.\\
\nonumber T(G_1 \oplus G_2,x,y)&=&T(G_1,x,y)\cdot T(G_2,x,y)\\
\nonumber T(E_1)&=&1;\\
T(\emptyset)&=&1; \label{rec_tutte}
\end{eqnarray}
However, we shall use in this paper the version of the Tutte
polynomial used by A.Sokal \cite{ar:Sokal2005a}:
\begin{equation}
Z(G,q,v)=\sum_{A\subseteq E}q^{k(A)}v^{|A|}
\end{equation}
The bivariate Sokal polynomial is co-reducible to the Tutte
polynomial via
\begin{equation}
T(G,x,y)=(x-1)^{-k(E)}(y-1)^{-|V|}Z(G,(x-1)(y-1),y-1)
\end{equation}
and has much recurrence relation
which does not distinguish whether the edge $e$ is a loop,
a bridge, or none of the two:
\begin{eqnarray}
\nonumber Z(G,q,v)&=&v \cdot Z(G_{/e},q,v) + Z(G_{-e},q,v)\\
\nonumber Z(G_1 \oplus G_2,q,v)&=&Z(G_1,q,v)\cdot Z(G_2,q,v)\\
\nonumber Z(E_1)&=&q;\\
Z(\emptyset)&=&1; \label{rec_sokal}
\end{eqnarray}

\paragraph{Bivariate chromatic polynomial}
K.Dohmen, A.P\"{o}nitz and P.Tittman in \cite{ar:DPT03} introduced
a polynomial $P(G,x,y)$
by splitting the available colors into colors for proper and colors
for arbitrary colorings.

We prove in this paper that this polynomial satisfies the following
recurrence relation:
\begin{eqnarray}\label{rec_DPT}
\nonumber P(G,x,y)&=&P(G_{-e},x,y) - P(G_{/e},x,y) + (x-y)\cdot P(G_{\dagger e},x,y)\\
\nonumber P(G_1 \oplus G_2,x,y)&=&P(G_1,x,y)\cdot P(G_2,x,y)\\
\nonumber P(E_1)&=&x;\\
P(\emptyset)&=&1;
\end{eqnarray}

\subsection{A most general edge elimination polynomial}
We define the most general confluent linear recurrence relation\footnote{
The first paper to study general conditions
under which linear recurrence relations define a graph invariant
is D.N. Yetter \cite{ar:Yetter90}.
},
which can be obtained on unlabeled graphs by introducing new
variables, and which does not distinguish between local properties
of the edge $e$ which is to be eliminated\footnote{
It is conceivable that recurrence relations with
various case distinctions depending on local
properties of $e$ and more variables give other ``most general'' polynomials.
This is the reason why we speak of ``a most general'' edge elimination
polynomial in the title of the paper.
}. 
We start with the recurrence relation
\begin{eqnarray}
\nonumber \xi(G)&=&w \cdot\xi(G_{-e}) + y\cdot \xi(G_{/e}) + z\cdot \xi(G_{\dagger e})\\
\nonumber \xi(G_1 \oplus G_2)&=&\xi(G_1)\cdot \xi(G_2)\\
\nonumber \xi(E_1)&=&x;\\
\xi(\emptyset)&=&1; \label{rec_ERP_with_w}
\end{eqnarray}
We prove:
\begin{thm}\label{thm_rec_ERP}
The recurrence relation (\ref{rec_ERP_with_w}) is confluent if and
only if one of the following conditions are satisfied:
\begin{eqnarray}
z = 0 \\
w = 1 \label{conf_condition2}
\end{eqnarray}
\end{thm}
Under the confluence condition (\ref{conf_condition2}), which allows
more general graph polynomial to be obtained, the recurrence
relation (\ref{rec_ERP_with_w}) is restricted to
\begin{eqnarray}
\nonumber \xi(G,x,y,z)&=&\xi(G_{-e},x,y,z) + y\cdot \xi(G_{/e},x,y,z) + z\cdot \xi(G_{\dagger e},x,y,z)\\
\nonumber \xi(G_1 \oplus G_2,x,y,z)&=&\xi(G_1,x,y,z)\cdot \xi(G_2,x,y,z)\\
\nonumber \xi(E_1,x,y,z)&=&x;\\
\xi(\emptyset,x,y,z)&=&1; \label{rec_ERP}
\end{eqnarray}

From this theorem one sees immediately that
the polynomial $\xi(G,x,y,z)$  gives, by choosing appropriate values
for the variables and simple prefactors,
the bivariate Sokal polynomial,
the bivariate matching polynomial and the bivariate
chromatic polynomial with all their respective substitution
instances, including the classical chromatic polynomial, the Tutte
polynomial, the vertex-cover and the independent set polynomial,
\cite{ar:DongHendyTeoLittle02,ar:GutmanHarary83}.
The latter two polynomials are already substitution instances
of the bivariate chromatic polynomial $P(G,x,y)$ of
\cite{ar:DPT03}.

In our next result
we give an explicit form of the polynomial
$\xi(G,x,y,z)$ using 3-partition expansion\footnote{
A more precise name would be ``Pair of two disjoint subsets expansion".
We chose the name 3-partition expansion, as any two disjoint subsets
induce a partition into three sets.
}:
\begin{thm}\label{thm_exp_ERP}
Let $G=(V,E)$ be a (multi)graph. Then the edge elimination
polynomial $\xi(G,x,y,z)$ can be calculated as
\begin{equation}\label{exp_ERP}
\xi(G,x,y,z) = \sum_{{\tiny
              (A \sqcup B) \subseteq E}}
              x^{k(A\sqcup B)-k_{cov}(B)}\cdot y^{|A|+|B|-k_{cov}(B)}\cdot
              z^{k_{cov}(B)}
\end{equation}
where by abuse of notation we use $(A\sqcup B)\subseteq E$ for
summation over subsets $A,B \subseteq E$, such that the subsets of
vertices $V(A)$ and $V(B)$, covered by respective subset of edges,
are disjoint: $V(A) \cap V(B) = \emptyset$; $k(A)$ denotes the
number of spanning connected components in $(V,A)$, and $k_{cov}(B)$
denotes the number of covered connected components, i.e. the
connected components of $(V(B),B)$.
\end{thm}
\begin{rem}
\label{msol}
From
Theorem \ref{thm_exp_ERP}  one can see that
$\xi(G,x,y,z)$
is a polynomial definable in
Monadic Second Order Logic, with quantification over
sets of edges ($MSOL_2$), where an order over vertices is to be used for
stating "number of connected sets", but the final result is
order-independent. We shall not use logic in the sequel of the paper.
For details the reader is referred to
\cite{ar:Makowsky01}.
\end{rem}

\subsection{Comparison with the weighted graph polynomial}
\label{se:noble}
The weighted graph polynomial
$U(G,\bar{x},y)$
introduced by
S.D. Noble and D.J.A. Welsh in \cite{ar:NobleWelsh99}
is defined for a graph $G=(V,E)$ as
$$
U(G,\bar{x},y)  =
\sum_{A \subseteq E} \prod_{i=1}^{|V|} x_i^{s(i,A)}  y^{|A|-r(A)}
$$
where $s(i,A)$ denotes the number of connected components
of size $i$ in the spanning subgraph $(V,A)$,
and $r(A)= |V|- k(A)$ is the rank of $(V,A)$.

The main difference between
$U(G,\bar{x},y)$  and
$\xi(G,x,y,z)$ is the number of variables,
which grows in the case of $U$ and is fixed in the case of $\xi$.
Furthermore, in the definition of $s(i,A)$ the numeric value of the
index of the variable $x_i$ is used.
This has as a consequence that one cannot freely rename the
variables of $U$. In $\xi$, as well as in all
graph polynomials definable in $MSOL_2$ in an order invariant way,
the variables can be renamed.
This allows one to show that
$U(G,\bar{x},y)$ is not an $MSOL_2$-definable polynomial.

$U(G,\bar{x},y)$
also gives
the Tutte polynomial and the matching polynomial as its substitution instances.
One can see that the
polynomial $U(G,\bar{x},y)$ distinguishes between graphs
for which $\xi(G,x,y,z)$ gives the same value.
As an example we look at the
trees shown on Fig. 1.
We do not know whether $\xi(G,x,y,z)$ can be obtained as a
substitution instance of $U(G,\bar{x},y)$.
\begin{center}
\includegraphics[width=2.2in]{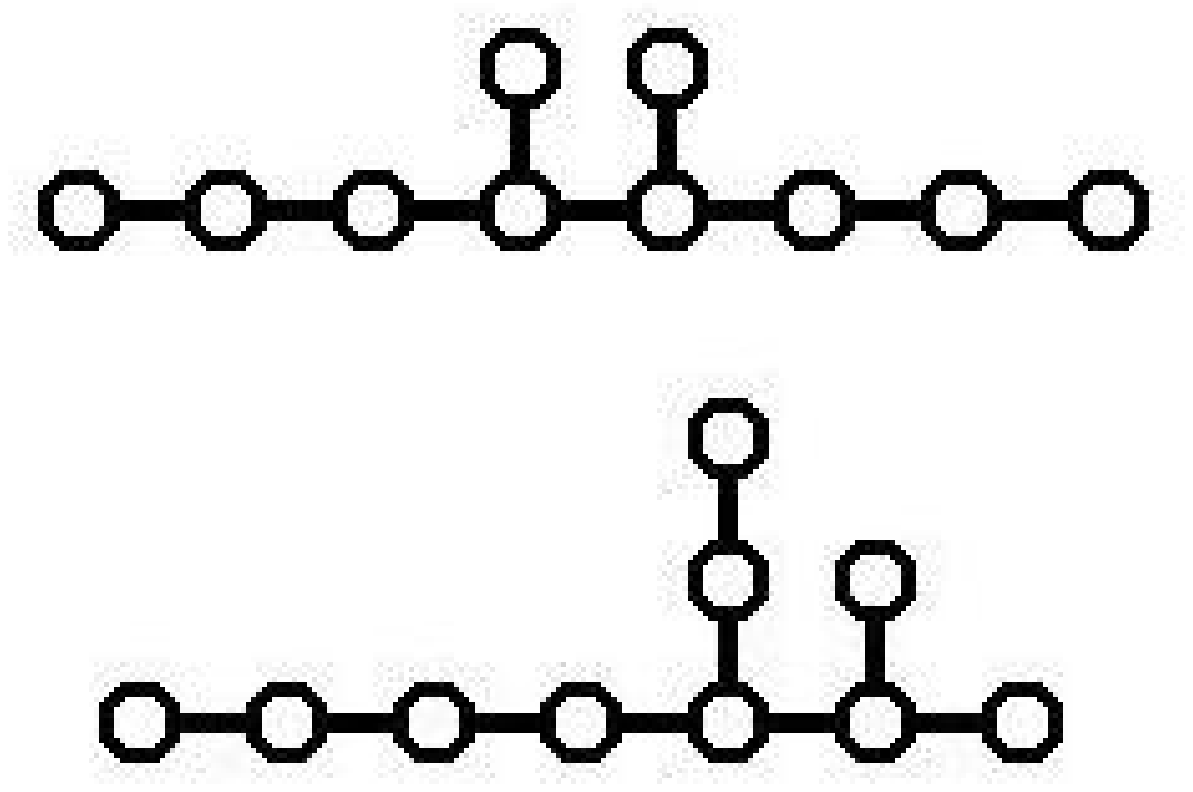}\\
\textbf{Fig. 1}: Non-isomorphic trees having the same
$\xi(G,x,y,z)$.
\end{center}

\subsection{A labeled version of $\xi$}
For edge-labeled\footnote{ In \cite{ar:BollobasRiordan99} they speak
of edge-colorings rather than edge-labelings. As we also discuss
chromatic polynomials we prefer our terminology as it avoids
confusions. } graphs we define the labeled version of our
polynomial: Let $G=(V,E,c)$ be an edge-labeled multigraph s.t.
$c:E\mapsto \Lambda$,
where $\Lambda$ is a set of labels, without any algebraic structure
defined over it,
and suppose that for each $\lambda \in \Lambda$ three elements
$w_{\lambda}$, $y_{\lambda}$ and $z_{\lambda}$ of a field are
chosen. Then using the same approach as for the unlabeled version,
we define a linear recurrence relation:
\begin{eqnarray}\label{col_rec_ERP_with_w}
\nonumber \xi_{lab}(G)&=&w_{c(e)} \cdot \xi_{lab}(G_{-e}) +
          y_{c(e)} \cdot \xi_{lab}(G_{/e}) +
          z_{c(e)} \cdot \xi_{lab}(G_{\dagger e})\\
\nonumber \xi_{lab}(G_1 \oplus G_2)&=&\xi_{lab}(G_1)\cdot \xi_{lab}(G_2)\\
\nonumber \xi_{lab}(E_1)&=&x;\\
\xi_{lab}(\emptyset)&=&1;
\end{eqnarray}
Note that we do not introduce weights on the vertices, as this would
make the definition of the edge contraction unclear,
unless we define an algebraic structure over $\Lambda$. For example,
if $\Lambda$ was a ring, we could define the label of the vertex
produced by a contraction of an edge $\{u,v\}$ to be the sum of the
labels of $u$ and $v$. In that case, we would get a generalization
of the weighted graph polynomial for labeled graphs
$W(G,\bar{x},y)$, also introduced by S.D.Noble and D.J.A.Welsh in
\cite{ar:NobleWelsh99}. However, this polynomial is "too strong", in
sense that it has the same definability problems as $U(G, \bar{x},y)$
discussed in Section \ref{se:noble}.

For the labeled case we prove 
\begin{thm}
\label{thm_lab_conf}
Every one of the conditions
\begin{gather}
\forall {e\in E} \left(z_{c(e)} = 0\right) \label{col_conf_condition1}\\
\forall {e\in E} \left(w_{c(e)} = 1\right) \wedge
\forall {e_1,e_2 \in E}\left(y_{c(e_1)}z_{c(e_2)} =
y_{c(e_2)}z_{c(e_1)}\right)\label{col_conf_condition2}
\end{gather}
is sufficient for the recurrence relation (\ref{col_rec_ERP_with_w})
being confluent.
\end{thm}
\begin{rem}
The conditions in Theorem
\ref{thm_lab_conf} are not necessary.
To see this we look at a graph with two connected
components and use condition
(\ref{col_conf_condition1}) for edges in the first component
and condition
(\ref{col_conf_condition2}) for edges in the second component.
\end{rem}

Under the confluence condition (\ref{col_conf_condition2}), which
allows more general graph polynomial to be obtained, the recurrence
relation (\ref{col_rec_ERP_with_w}) is restricted to
\begin{eqnarray}\label{col_rec_ERP}
\nonumber \xi_{lab}(G)&=&\xi_{lab}(G_{-e}) +
          y \cdot t_{c(e)} \cdot \xi_{lab}(G_{/e}) +
          z \cdot t_{c(e)} \cdot \xi_{lab}(G_{\dagger e})\\
\nonumber \xi_{lab}(G_1 \oplus G_2)&=&\xi_{lab}(G_1)\cdot \xi_{lab}(G_2)\\
\nonumber \xi_{lab}(E_1)&=&x;\\
\xi_{lab}(\emptyset)&=&1;
\end{eqnarray}
where
\begin{eqnarray}\label{eqn_t_y_z}
\nonumber y_{c(e)} = y\cdot t_{c(e)}\\
z_{c(e)} = z\cdot t_{c(e)}
\end{eqnarray}
$x$, $y$ and $z$ are unlabeled variables, and $\bar{t}$ is the
unique solution of (\ref{eqn_t_y_z}). Like the unlabeled case, we
also introduce the explicit form:
\begin{thm}
The expression
\begin{equation}\label{col_exp_ERP}
\xi_{lab}(G,x,y,z,\bar{t}) = \sum_{{\tiny
              (A \sqcup B) \subseteq E}}
              x^{k(A\sqcup B)}\left(\prod_{e\in A \sqcup B}(y t_{c(e)})\right)
              \left(\frac{z}{xy}\right)^{k_{cov}(B)}
\end{equation}
defines the same graph polynomial as the recurrence relation
(\ref{col_rec_ERP}).
\end{thm}
Note that the degree of $x$ and $y$ in the denominator does never
exceed the degree of the respective variable in the nominator.

\begin{rem}
The \emph{labeled Sokal polynomial}
\cite{ar:Sokal2005a}, Zaslavsky's \emph{normal function of the
colored matroid} \cite{ar:Zaslavsky92}, Heilmann and Lieb's
\emph{labeled matching polynomial} \cite{ar:HeilmannLieb72},
and the chain polynomial \cite{ar:ReadWhitehead99,ar:Traldi02} are
substitution instances of $\xi_{lab}(G,x,y,z,\bar{t})$
up to a simple prefactor.
\end{rem}

The remainder of the paper is organized as follows: in Section
\ref{sec_rec_DPT} we prove the recurrence relation of the
generalized chromatic polynomial. In Section
\ref{sec_gen_rec_ERP} we establish the most general linear
recurrence relation with respect to the three edge elimination
operations, restrict it to be multiplicative and confluent, and then
prove the confluence property of the resulting function. In
Section \ref{sec_gen_exp_ERP} we establish the explicit function as
in Theorem \ref{thm_exp_ERP}, and prove that it defines the same
polynomial. The section \ref{sec_col_ERP} expands our results to the
edge-labeled graphs. The section \ref{sec_subst_ERP} contains
examples of known graph polynomials which can be obtained as
substitution instances of the edge elimination polynomial. Finally,
in Section \ref{sec_comp_ERP} we deal with the complexity of its
computation.

\section{The recursive definition of the \\ generalized chromatic
polynomial}\label{sec_rec_DPT}

Recall the definition given by K.Dohmen, A.P\"{o}nitz and P.Tittman
in \cite{ar:DPT03}: There are two disjoint sets of colors $Y$ and
$Z$; a generalized coloring of a graph $G=(V,E)$ is a map
$\phi:V\mapsto (Y \sqcup Z)$ such that for all $\{u,v\}\in E$, if
$\phi(u) \in Y$ and $\phi(v) \in Y$, then $\phi(u) \neq \phi(v)$
(The set $Y$ is called therefore "the proper colors"). For two
positive integers $x>y$, the value of the polynomial is the number
of generalized colorings of $G$ by $x$ colors, $y$ of them are
proper. We enhance this definition to multigraphs by the following:
\begin{enumerate}
\item A self-loop can be colored only by a color in $X \setminus Y$;
\item A multiple edge does not affect colorings.
\end{enumerate}
Let $G=(V,E)$ be a graph, and $P(G,x,y)$ be the number of
generalized colorings defined above. Let $v\in V$ be any vertex. We
denote by $P^v(G,x,y)$ the number of generalized colorings of $G$,
when $v$ is not colored by a proper color, i.e. $\phi(v)\in
X\setminus Y$.
\begin{prop}\label{prop_colv}
$P^v(G,x,y) = (x-y)\cdot P(G_{-v},x,y)$, where $G_{-v}$ denotes the
subgraph of $G$ induced by $V\setminus \{v\}$.
\end{prop}
\begin{proof}
By inspection: the vertex $v$ can have any color in $X\setminus Y$,
and the coloring of the remainder does not depend on it. \end{proof}
Let $e=\{u,v\} \in E$ be any edge of $G$, which is not a self-loop
and not a multiple edge. Consider the number of colorings of
$G_{-e}$. Any such coloring is either a coloring of $G$, or a
coloring of $G_{/e}$, when the vertex $u=v$, which is produced by
the contraction, is colored by a proper color. Together with
Proposition \ref{prop_colv}, that raises:
\begin{equation}
P(G,x,y)=P(G_{-e},x,y) - P(G_{/e},x,y) + (x-y)\cdot P(G_{\dagger
e},x,y)
\end{equation}
One can easily check that this equation is satisfied also for loops
and multiple edges. Together with the fact that a singleton can be
colored by any color, and the fact that the number of colorings is
multiplicative, this proves the recursive definition
(\ref{rec_DPT}).

\section{The most general recurrence relation}
\label{sec_gen_rec_ERP}
We are looking for the most general linear recurrence relation with
respect to edge deletion, edge contraction and edge extraction operation
that can be obtained by introducing new variables. Recall that we
are interested in a \emph{graph invariant}, e.i. the resulting
function should not depend on the order of graph deconstruction.
Moreover, this invariant should be a multiplicative graph
polynomial.

From this consideration alone we obtain the initial condition and
the product rule:
\begin{eqnarray}
\nonumber \xi(G_1 \oplus G_2)&=&\xi(G_1)\cdot \xi(G_2)\\
\xi(\emptyset)&=&1; \label{rec_ERP_init}
\end{eqnarray}
Indeed, the disjoint union with an empty set gives the same graph,
so the resulting function should also remain the same.\\

At this stage, we formulate the edge elimination rule introducing a
new variable wherever we can. We set
\begin{eqnarray}
\nonumber \xi(G,x,y,z,t)&=&t \cdot \xi(G_{-e},x,y,z,t) +
y\cdot \xi(G_{/e},x,y,z,t) + z\cdot \xi(G_{\dagger e},x,y,z,t)\\
\nonumber \xi(E_1,x,y,z,t)&=&x;\\
\end{eqnarray}

Let $G$ be a graph as presented on Fig. 2. Note that the subgraphs
$H_1$, $H_{1-u}$, $H_2$ and $H_{2-w}$ can be different and have (in
general) different $\xi$.
\begin{center}
\includegraphics[width=2.5in]{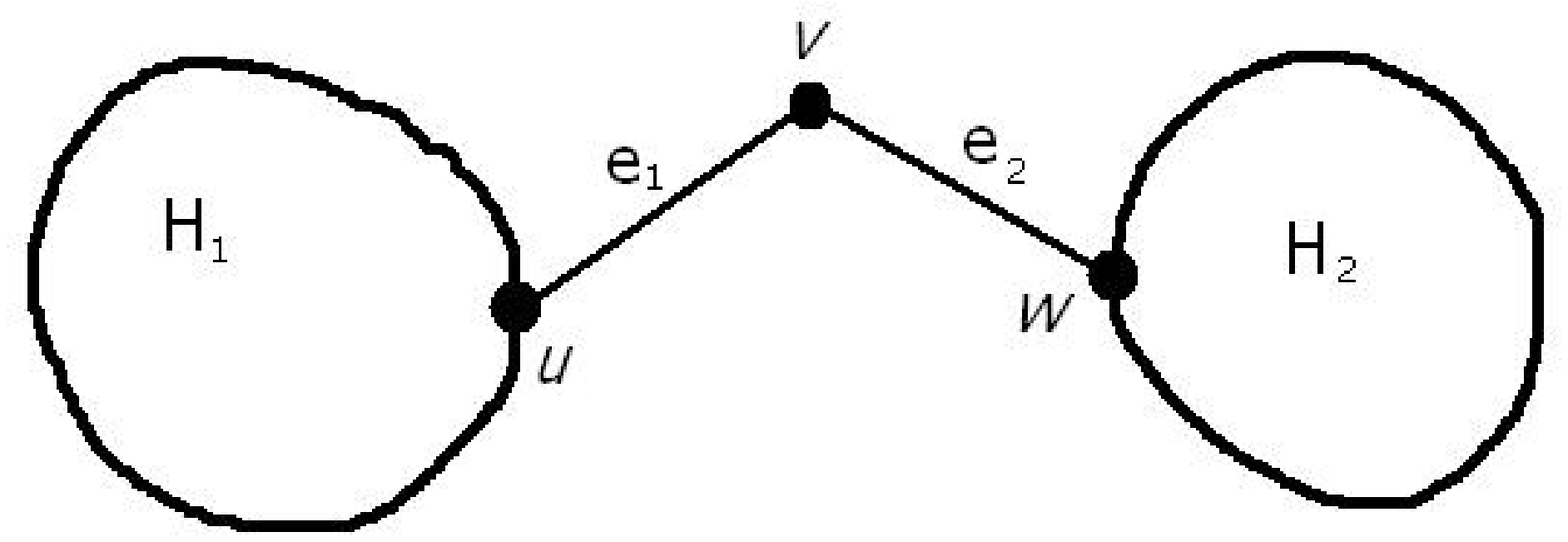}\\
\textbf{Fig. 2}: The graph for testing of the confluence property.
\end{center}
Since we are looking for a graph invariant, we must obtain the same
result by applying the edge elimination rule first on the edge $e_1$
and then on the edge $e_2$, as in case when we apply the edge
elimination rule first on the edge $e_2$ and then on the edge $e_1$.
\begin{eqnarray}
\nonumber \xi(G) &=& t \cdot \xi(G_{-e_1}) + y\cdot \xi(G_{/e_1}) + z\cdot \xi(G_{\dagger e_1}) =\\
\nonumber &=& t \cdot \xi(H_1) \cdot \left[ x \cdot t \cdot \xi(H_2)
+ y \cdot \xi(H_2) + z \cdot \xi(H_{2-w}) \right] + \\
\nonumber & & y \cdot \left[ t \cdot \xi(H_1)\xi(H_2) + y \cdot
\xi(G_{/e_1/e_2}) + z\cdot \xi(H_{1-u})\xi(H_{2-w})\right] + \\
& & z \cdot \xi(H_{1-u})\xi(H_2)
\end{eqnarray}
On the other hand,
\begin{eqnarray}
\nonumber \xi(G) &=& t \cdot \xi(G_{-e_2}) + y\cdot \xi(G_{/e_2}) + z\cdot \xi(G_{\dagger e_2}) =\\
\nonumber &=& t \cdot \xi(H_2) \cdot \left[ x \cdot t \cdot \xi(H_1)
+ y \cdot \xi(H_1) + z \cdot \xi(H_{1-u}) \right] + \\
\nonumber & & y \cdot \left[ t \cdot \xi(H_1)\xi(H_2) + y \cdot
\xi(G_{/e_1/e_2}) + z\cdot \xi(H_{1-u})\xi(H_{2-w})\right] + \\
& & z \cdot \xi(H_{2-w})\xi(H_1)
\end{eqnarray}
Hence, we have either $z=0$ or $t=1$ as a necessary condition of the
confluence. In case of $z=0$ the resulting function is a
substitution instance of the Sokal polynomial:
\begin{equation}
\xi(G,x,y,0,t)=t^{|E|}\cdot Z(G,x,\frac{y}{t})
\end{equation}
We leave the proof of that fact to the reader. Since the Sokal
polynomial can be also obtained when $t=1$, the latter case is
considered more general, and it will be further investigated. That
brings us back to the recurrence relation (\ref{rec_ERP}). To
complete the proof of Theorem \ref{thm_rec_ERP}, we need now to show
that
\emph{the recurrence relation (\ref{rec_ERP}) is confluent}.\\

It is enough to prove that any two steps of the graph decomposition
using (\ref{rec_ERP}) are interchangeable. This includes two parts:
\begin{itemize}
\item Decomposition of a graph by elimination of any two edges in different order;
\item Edge elimination and disjoint union.
\end{itemize}
The proof of both parts is rather technical and left to the reader.

\section{The explicit form or the polynomial $\xi(G,x,y,z)$}
\label{sec_gen_exp_ERP}
In this section we prove Theorem \ref{thm_exp_ERP}. In order to do
so, we need to show that
\begin{itemize}
\item The expression (\ref{exp_ERP}) satisfies the initial
conditions of (\ref{rec_ERP});
\item The expression (\ref{exp_ERP}) is multiplicative;
\item The expression (\ref{exp_ERP}) satisfies the edge elimination
rule of (\ref{rec_ERP}).
\end{itemize}
Then by induction on the number of edges in $G$ the theorem holds.
The first fact is trivial; the second one can be easily checked by
reader. Indeed, the summation over subsets of edges of $G(V,E) =
G_1(V_1,E_1) \oplus G_2(V_2,E_2)$ can be regarded as a summation
over the subsets of $E_1$, and then independently over the subsets
of $E_2$. Therefore, we just need to prove that
\begin{lem}
The explicit expression given by (\ref{exp_ERP}) satisfies the edge
elimination rule of (\ref{rec_ERP}).
\end{lem}
\begin{proof}
Let $G=(V,E)$ be the (multi)graph of interest. Let $N(G)$ be defined
as
\begin{equation}
N(G,x,y,z) = \sum_{{\tiny
              (A \sqcup B) \subseteq E}}
              x^{k(A\sqcup B)-k_{cov}(B)}\cdot y^{|A|+|B|-k_{cov}(B)}\cdot
              z^{k_{cov}(B)}
\end{equation}
where $k(A)$ denotes the number of connected components in $(V,A)$,
and $k_{cov}(B)$ denotes the number of the connected components of
$(V(B),B)$, where $V(B)\subseteq V$ are the vertices covered by the
edges of $B$. Let $e$ be the edge we have chosen to reduce. Any
particular choice of $A$ and $B$ can be regarded as a
vertex-disjoint edge coloring in 2 colors A and B, when part of the
edges remains uncolored. We divide all the coloring into three
disjoint cases:
\begin{itemize}
\item Case 1: $e$ is uncolored;
\item Case 2: $e$ is colored by $B$, and it is the only edge of a colored
connected component;
\item Case 3: All the rest. That means, $e$ is colored by $A$, or $e$ is
colored by $B$ but it is not the only edge of a colored connected
component.
\end{itemize}
In the case 1, we just sum over colorings of $G_{-e}$:
\begin{equation}
N_1(G)=\sum_{(A \sqcup B) \models ~Case~ 1}
x^{k(A\sqcup B)-k_{cov}(B)}\cdot y^{|A|+|B|-k_{cov}(B)}\cdot
z^{k_{cov}(B)}
=N(G_{-e})
\end{equation}
In the case 2, the edge $e$ is a connected component of $(V(B),B)$.
Therefore, if we analyze now $N(G_{\dagger e})$, we will get
\begin{itemize}
\item The number of edges colored by $A$ is the same;
\item The number of edges colored by $B$ is reduced by one;
\item The total number of colored connected components is reduced by one;
\item The number of covered connected components colored $B$ is
reduced by one;
\end{itemize}
This gives us
\begin{equation}
N_2(G)=\sum_{(A \sqcup B) \models ~Case~ 2}
x^{k(A\sqcup B)-k_{cov}(B)}\cdot y^{|A|+|B|-k_{cov}(B)}\cdot
z^{k_{cov}(B)}
=z \cdot N(G_{\dagger e})
\end{equation}
And finally, in the case 3, $e$ is a part of a bigger colored
connected component, or it is alone a connected component colored by
$A$. In this case, we analyze the colorings of $G_{/e}$:
\begin{itemize}
\item Either $|A|$ or $|B|$ is reduced by 1, the other remained the
same;
\item The total number of colored connected components
remained the same;
\item The number of covered connected components colored $B$
remained the same.
\end{itemize}
According to the above,
\begin{equation}
N_3(G)=\sum_{(A \sqcup B) \models ~Case~ 3}
x^{k(A\sqcup B)-k_{cov}(B)}\cdot y^{|A|+|B|-k_{cov}(B)}\cdot
z^{k_{cov}(B)}
=y \cdot N(G_{/e})
\end{equation}
which together with $N(G) = N_1(G) + N_2(G) + N_3(G)$ completes the
proof.
\end{proof}
\section{The edge elimination polynomial of a labeled graph}
\label{sec_col_ERP}
To obtain the edge-labeled version of our polynomial, we use the
same approach as in Section \ref{sec_gen_rec_ERP}: we are looking
for a multiplicative graph invariant satisfying linear recurrence
relation with respect to the edge elimination operations. We start
with
\begin{eqnarray}
\nonumber \xi_{lab}(G_1 \oplus G_2)&=&\xi_{lab}(G_1)\cdot \xi_{lab}(G_2)\\
\xi_{lab}(\emptyset)&=&1; \label{rec_col_ERP_init}
\end{eqnarray}
and define an edge elimination rule introducing a new variable
wherever we can (now every variable has index $e$ for the edge which
is currently being eliminated\footnote{We do not use an index for
vertices, because the vertex set of the graph is being changed
during decomposition. For this reason, we cannot call this
polynomial "the most general"}.
\begin{eqnarray}\label{rec_col_ERP_step}
\nonumber \xi_{lab}(G)&=&w_e \cdot \xi(G_{-e}) + y_e \cdot \xi_{lab}(G_{/e}) + z_e \cdot \xi_{lab}(G_{\dagger e})\\
\nonumber \xi_{lab}(E_1)&=&x;\\
\end{eqnarray}
The same considerations as in Section \ref{sec_gen_rec_ERP}, using
the same graph (Fig. 1), we get that the recursion
(\ref{rec_col_ERP_step}) is confluent when either $z_e = 0$ or $w_e
= 1$ and $y_{e_1}z_{e_2} = y_{e_2}z_{e_1}$. One can expand this
result to any two edges of a connected component. Since the graph in
general can be connected, and our recurrence relation should be
confluent for \emph{every} graph, we get the following restrictions:
\begin{itemize}
\item $z_e = 0$ for \emph{every} edge $e$, or
\item $w_e = 1$ , $y_e = y \cdot t_e$ and
$z_e = z \cdot t_e$ for \emph{every} edge $e$ (here $y$ and $z$ do
not depend on $e$).
\end{itemize}
In the first case we obtain an instance of the labeled Sokal
polynomial:
\begin{equation}
\xi_{lab}(G,\bar{w},x,\bar{y},\bar{0})=\left(\prod_{e\in E} w_e
\right)\cdot Z(G,q,\bar{v})
\end{equation}
where $q=x$ and $v_e=\frac{y_e}{w_e}$. Since the Sokal polynomial
can be also obtained when $\bar{w}=\bar{1}$, the latter case is
considered more general. That brings us to the recurrence relation
(\ref{col_rec_ERP}). We have now to prove two propositions:
\begin{prop}
The recurrence relation (\ref{col_rec_ERP}) is confluent.
\end{prop}
\begin{prop}
The formula (\ref{col_exp_ERP}) defines the same polynomial as the
recurrence relation (\ref{col_rec_ERP}).
\end{prop}
Both the proofs are similar to the respective unlabeled version and
left to the reader.
\section{Application to some known graph polynomials}
\label{sec_subst_ERP}
In this section we present different known graph polynomials as
substitution instances of $\xi(G)$ and $\xi_{lab}(G)$. Two issues
should be addressed here:
\subparagraph{Zero coefficients:}%
When some of the arguments $x$, $y$ or $z$ of our polynomial is
zero, we generally get 0 in all the summands that contain this
variable in some positive power, and an uncertainty of kind $0^0$ in
all the summands that contain it in power 0. However, as of being a
polynomial, our function is continuous, and thus we can use the fact
that for any nonnegative integer $k$,
$$x^k|_{x=0} = \lim_{x\rightarrow 0} x^k = \left\{
              \begin{array}{ll}
              1~~~if~k=0 \\
              0~~~otherwise
              \end{array}\right.$$
Hence, if in our substitution some variable turns 0, the value of
the resulting polynomial is still well-defined.
\subparagraph{Multiple edges and loops:}%
Some of the graph polynomials are defined only for simple and
loop-free graphs. However, their definition can be easily
generalized to multigraphs, such that the equality holds in case of
a simple input graph.

\subsubsection*{Labeled versions of the Tutte polynomial}
\begin{prop}
The Sokal polynomial (in both unlabeled and labeled versions) can be
obtained by
$$Z(G,q,v)=\xi(G,q,v,0)$$
$$Z_{lab}(G,q,\bar{v}) = \xi_{lab}\left(G,q,1,0,\bar{v}\right)$$
in particular, the chromatic polynomial can be obtained by
$$\chi(G,\lambda) = Z(G,\lambda,-1) = \xi(G,\lambda,-1,0)$$
\end{prop}
\begin{proof} 
By inspection of summands with $B=\emptyset$. All the
other summands are eliminated by $z=0$.\end{proof}%
By a simple substitution of variables, we get the following three
corollaries:
\begin{cor}
The classical Tutte polynomial can be obtained by
$$T(G,x,y)=(x-1)^{-k(E)}\cdot(y-1)^{-|V|}\cdot\xi\left(G,(x-1)(y-1),(y-1),0\right)$$
\end{cor}
Recall that $r(S) = |V| -k(S)$ is the rank
of the spanning subgraph with edge set $S$.
The Zaslavsky's normal function of the colored matroid, applied to a
graph $G=(V,E)$ with edge coloring function 
$c:E\mapsto \Lambda$,
is defined by
$$
R(G,c)=\sum_{S\subseteq E}
\left(\prod_{e\in S} x_{c(e)}\right)
\left(\prod_{e\not\in S}y_{c(e)}\right)
(x-1)^{r(E)-r(S)}
(y-1)^{|S|-r(S)}.
$$
\begin{cor}
The Zaslavsky's normal function of the edge-colored graph
can be obtained by %
\begin{eqnarray}
\nonumber
R(G,c)=
(x')^{-k(E)} \cdot (y')^{-|V|} \cdot%
\left(
\prod_{e\in E}y_{c(e)}
\right) \cdot%
\xi_{lab}\left(G,x'y',y',0,\bar{t}\right)
\end{eqnarray}
where $t_{c(e)}=\frac{x_{c(e)}}{y_{c(e)}}$
and $x'=x-1, y'=y-1$.
\end{cor}

The chain polynomial $Ch(G,\omega, \bar{u})$ was first introduced
in \cite{ar:ReadWhitehead99} and
can also be defined, cf. \cite{ar:Traldi02},
as 
$$
Ch(G, \omega, \bar{u}) =
\sum_{S \subseteq E} (1- \omega)^{|S|- r(S)} \prod_{e \in E-S} u_e.
$$
From \cite{ar:Traldi02} we get 
\begin{cor}
The chain polynomial can be obtained by
$$
Ch(G,\omega, \bar{ua)} =
\left(
\prod_{e \in E} u_e
\right)
\cdot
(1- \omega)^{-|V|} \cdot \xi_{lab}(G, 1 - \omega, 1, 0, \bar{v})
$$
where $v_e = \frac{1 - \omega}{u_e}$.
\end{cor}

\subsubsection*{Matching polynomials}
The next two propositions deal with various forms of matching
polynomials:
\begin{prop}
The generalized matching polynomial (\ref{general_matching_eqn}) can
be obtained by
$$M(G,x,y)=\sum_{i=0}^n a_i x^{n-2i} y^i = \xi(G,x,0,y)$$
In particular, the generating matching polynomial is
$g(G,x)=\xi(G,1,0,x)$ and the defect matching polynomial is
$\mu(G,x)=\xi(G,x,0,-1)$
\end{prop}
\begin{prop} The original Heilmann and Lieb's multivariate matching
polynomial introduced in \cite{ar:HeilmannLieb72} can be obtained by
\begin{equation}
\label{HL} M_{col}(G, \bar{x}, \bar{y}) =\sum_{{\tiny
\begin{array}{cc}
              M\subseteq E, \\
              M~is~a~matching
              \end{array}}}
\prod_{e=\{u,v\}\in M}y_e x_u x_v = \xi_{lab}(G,1,0,1,\bar{t})
\end{equation}
where $t_{e}=y_e x_u x_v$ for every edge $e=\{u,v\}$.
\end{prop}
\begin{proof}
By inspection of non-zero summands of $\xi(G)$. There should be no
edges in $A$, and every edge of $B$ should be in different connected
component, so $B$ has to be a matching.
\end{proof}
Finally, the Dohmen-P\"{o}nitz-Tittman generalization of the
chromatic polynomial \cite{ar:DPT03} is also a substitution instance
of $\xi(G)$:
\begin{prop}
$$P(G,x,y)=\xi(G,x,-1,x-y)$$
\end{prop}
\begin{proof}
Using recursion scheme (\ref{rec_DPT}), by induction on number of
edges $|E|$.
\end{proof}
Using results of K.Dohmen, A.P\"{o}nitz and P.Tittman
\cite{ar:DPT03} we can also derive that the independence polynomial
(which is a substitution instance of $P(G,x,y)$) is also a
substitution instance of $\xi(G)$.
\section{Computational complexity of $\xi(G)$}\label{sec_comp_ERP}

In this section we analyze the complexity of computation of $\xi(G)$
and $\xi_{lab}(G)$. In general, these polynomials are
$\sharp\mathbf{P}$-hard to compute, as every instance stated in the
previous section is $\sharp\mathbf{P}$-hard.
Recall that, according to Remark \ref{msol},
the formulas (\ref{exp_ERP}) and (\ref{col_exp_ERP}) can be
used to give an order invariant definition in
Monadic Second Order Logic, with quantification over
sets of edges, and an auxiliary order.

Hence, due to the general theorem from
\cite{ar:Makowsky01,ar:MakowskyTARSKI}, we have
\begin{prop}
\label{prop:tw}
$\xi(G)$ and $\xi_{lab}(G)$ are polynomial time computable on graphs
of tree-width at most $k$ where the exponent of the run time is
independent of $k$.
\end{prop}
Recall also from Remark \ref{msol} that
the weighted graph polynomial
$U(G,\bar{x},y)$ is not definable using
$MSOL$, and, hence, the results of
\cite{ar:Makowsky01,ar:MakowskyTARSKI}
are not applicable.
Indeed, the run time of the algorithm introduced by C.Noble in
\cite{ar:Noble06} for graphs of tree width at most $k$ is
polynomial, but its highest degree depends on $k$.

The drawback of the general method of
\cite{ar:Makowsky01,ar:MakowskyTARSKI}
lies in the huge hidden
constants, which make it practically unusable.
However, an explicit dynamic algorithm for computing the polynomial
$\xi_{lab}(G)$ on graphs of bounded tree-width, given the tree
decomposition of the graph, where the constants are simply
exponential in $k$, can be constructed along the same ideas as
presented in \cite{ar:Traldi03,ar:FischerMakowskyRavve2006}.

\section{Open questions}
\paragraph{Difficult point property:} In general, the computation of
$\xi(G)$ is $\sharp\mathbf{P}$-hard. However, we know that for some
$x$, $y$ and $z$, it can be easy. F.Jaeger, D.Vertigan and D.Welsh
define in \cite{ar:JaegerVertiganWelsh90} the set of the points of the
$(x,y)$-plane in which the computation of the induced graph invariants
of the bivariate Tutte
polynomial is easy, proving that the remaining points are
$\sharp\mathbf{P}$-hard. Similar theorems have been proven for the
interlace polynomial, the cover polynomial and the colored Tutte
polynomial
\cite{ar:BlaeserDell07,ar:BlaeserHoffmann07,ar:BlaeserDellMakowsky07a}.
In all the cases the ``easy" points lay in a
semi-algebraic subset of the polynomial domain of lower dimension.
\begin{question}
Describe the set of points for which the induced graph
invariants $\xi(G)$ and $\xi_lab(G)$ are easy to compute.
\end{question}

\paragraph{Distinctive power:} We know that the polynomial $\xi(G)$
has at least the same distinctive power as the Tutte polynomial and
the bivariate chromatic polynomial together, but more than every one
of them individually. Indeed, since $T(G,x,y)$ and $P(G,x,y)$ are both
substitution instances of $\xi(G)$, if $\xi(G)$ coincides for two
graphs, so do $T(G,x,y)$ and $P(G,x,y)$. On the other hand, we do
not know whether $\xi(G)$ has more distinctive power.
\begin{question}
Are there two graphs $G_1, G_2$ such that for all $x,y$ we have
$$
T(G_1,x,y) = T(G_2,x,y) \mbox{  and  }  P(G_1,x,y) = P(G_2,x,y)
$$
but such that for some $x,y,z~$
$$\xi(G,x,y,z) \neq \xi(G_2,x,y,z)?$$
\end{question}

\paragraph{Complexity on graphs of bounded clique-width:}
We have seen in Proposition \ref{prop:tw}
that for graphs of tree-width at most $k$ computing the
edge reduction polynomials $\xi(G)$ and $\xi_{lab}(G)$
is fixed parameter tractable (FPT) in the sense of
\cite{bk:DowneyFellows99,bk:FlumGrohe2006}.
Another graph parameter, introduced in \cite{ar:CourcelleOlariu00}
and discussed there is the clique-width.
It is open
whether the Tutte polynomial is fixed parameter tractable
for graphs of clique-width at most $k$,
\cite{pr:GimenezHlinenyNoy2005,pr:MRAG06}.

\begin{question}
Are the polynomials
$\xi(G)$ and $\xi_{lab}(G)$ fixed parameter tractable
for graph  classes of bounded clique-width?
\end{question}

\subsection*{Acknowledgments}
We would like to thank B. Courcelle for his
comments on an early version of this paper.

\small

\end{document}